\documentclass[12pt,a4paper]{article}
\usepackage{latexsym,amssymb,amsfonts,amsmath,amsthm,nccmath,float,graphicx,multirow,enumitem,setspace,bm,authblk,tabularx,longtable,mathrsfs,needspace}
\usepackage{booktabs}
\usepackage[dvipsnames]{xcolor}
\usepackage[hidelinks,hyperfootnotes=false]{hyperref}
\usepackage[margin=1.83cm]{geometry}
\usepackage{placeins}
\hypersetup{
	colorlinks,
	linkcolor={red!80!black},
	citecolor={blue!80!black},
	urlcolor={blue!80!black}
}

\newtheorem{Theorem}{Theorem}[section]
\newtheorem{Proposition}{Proposition}[section]
\newtheorem{Corollary}{Corollary}[section]
\newtheorem{Conjecture}{Conjecture}[section]

\newtheorem{Lemma}{Lemma}[section]

\newcommand{\Z}{\mathbb{Z}}

\let\oldproofname=\proofname
\renewcommand{\proofname}{\rm\bf{\oldproofname}}

\allowdisplaybreaks

\begin{document}
	
\title{Existence of strong starters in finite abelian groups}

\author[a]{Xiaomiao Wang}
\author[b,c]{Tao Feng}
\author[b]{Hengrui Liu}
\author[b]{Weiyi Sun}
\affil[a]{School of Mathematics and Statistics, Ningbo University, Ningbo, 315211, P.R. China}
\affil[b]{School of Mathematics and Statistics, Beijing Jiaotong University, Beijing, 100044, P.R. China}
\affil[c]{Hebei Provincial Key Laboratory of Mathematical Theory and Analysis for Network and Data Science, Beijing Jiaotong University, Beijing, 100044, P.R. China}

\renewcommand*{\Affilfont}{\small\it}
\renewcommand\Authands{ and }

\affil[ ]{wangxiaomiao@nbu.edu.cn; tfeng@bjtu.edu.cn; henryleo@bjtu.edu.cn; vesunn@bjtu.edu.cn}
\date{}
\maketitle

\begin{abstract}
Let $(G,+)$ be a finite abelian group of odd order with identity $0_G$. A strong starter in $G$ is a set of unordered pairs $\{\{x_i,y_i\}:1\le i\le (|G|-1)/2\}$ such that these pairs partition $G\setminus\{0_G\}$, the set $\{\pm(x_i-y_i):1\le i\le (|G|-1)/2\}$ equals $G\setminus\{0_G\}$, and the sums $x_i+y_i$, $1\le i\le (|G|-1)/2$, are distinct and nonzero. Horton (1989) conjectured that every nontrivial finite abelian group of odd order admits a strong starter except for $\Z_3$, $\Z_5$, $\Z_9$, and $\Z_3\times\Z_3$. In this paper, we confirm Horton's conjecture. As applications, we determine which cycles of odd order or twice an odd order admit optimal extended irregular dominating sets. We also establish the existence of three pairwise orthogonal starters in finite abelian groups of odd order, and the existence of cyclic Room squares.
\end{abstract}


\footnotetext{Supported by Ningbo Natural Science Foundation under Grant 2024J018, and NSFC under Grant 12271023}

\section{Introduction}

Let $(G,+)$ be a finite abelian group of odd order with identity $0_G$. A \emph{starter} in $G$ is a set of unordered pairs $\mathcal S=\{\{x_i,y_i\}:1\le i\le (|G|-1)/2\}$ such that 
\begin{enumerate}
\item[(1)] $\{x_i:1\le i\le(|G|-1)/2\}\cup \{y_i:1\le i\le(|G|-1)/2\} = G\setminus \{0_G\}$;
\item[(2)] $\{\pm(x_i-y_i):1\le i\le(|G|-1)/2\} = G\setminus \{0_G\}$.
\end{enumerate}
A starter is called \emph{strong} if the sums $x_i+y_i$, for $1\le i\le (|G|-1)/2$, are pairwise distinct and nonzero. 

Strong starters were first introduced by Stanton and Mullin \cite{MS} to construct Room squares, which are named after T.G.~Room and have their origins in Kirkman's Fifteen Schoolgirls Problem (cf. \cite{DS1992}). Early constructions of strong starters focused primarily on finite fields. Mullin and Nemeth \cite{MN} gave a cyclotomic construction of a strong starter in the additive group of $\operatorname{GF}(q)$ whenever $q=2^k t+1$, where $k$ is a positive integer and $t>1$ is odd. Thus, their construction does not cover prime powers of the form $2^k+1$, which include $9$ and the Fermat primes. Horton \cite{H1981} and Dinitz \cite{D1979} generalized Mullin and Nemeth's cyclotomic construction, producing families of pairwise orthogonal strong starters in the additive groups of finite fields. These families were also used to construct pairwise orthogonal idempotent symmetric Latin squares. Chong and Chan \cite{ChongChan} provided a construction for strong starters in $\operatorname{GF}(q)$ for Fermat primes $q>5$. Taken together, these constructions show that the additive group of $\operatorname{GF}(q)$ admits a strong starter for every odd prime power $q$ except $q=3,5,9$.

Let $\Z_v$ denote the additive group of integers modulo $v$. Dinitz and Stinson \cite{DS80} used the hill-climbing algorithm described in  \cite{DS} to search for strong starters in cyclic groups, establishing their existence in
$\mathbb Z_v$ for every odd integer $v$ with $7\le v\le999$ and $v\ne9$. 

\begin{Lemma}\label{999}{\rm \cite{DS80}}
For every odd integer $v$ with $7\le v\le 999$ and $v\ne9$, there exists a strong starter in $\Z_v$.
\end{Lemma}

Gross \cite{G74} gave a multiplication construction for strong starters in direct sums of finite abelian groups. Gross and Leonard \cite{GL} subsequently established a recursive construction: if $H$ is a subgroup of an abelian group $G$, both $H$ and $G/H$ admit strong starters, and there is a permutation $\pi$ of $H$ such that the maps $x\mapsto\pi(x)+x$ and $x\mapsto\pi(x)-x$ are also
permutations of $H$, then $G$ admits a strong starter. Building on this recursive construction, Horton \cite[Theorem 2.12]{H} proved the following existence result.

\begin{Lemma}\label{known}{\rm \cite{H}}
Let $G$ be a finite abelian group of order $v>5$ with $\gcd(v,6)=1$. Then there exists a strong starter in $G$.
\end{Lemma}

Up to isomorphism, the finite abelian groups of orders $3$, $5$, and $9$ are $\Z_3$, $\Z_5$, $\Z_9$, and $\Z_3\times\Z_3$. A direct check shows that none of these groups admits a strong starter. Horton formulated the following conjecture.

\begin{Conjecture}\label{conj}{\rm \cite{H}}
Every nontrivial finite abelian group of odd order admits a strong starter, except for $\mathbb Z_3$, $\mathbb Z_5$, $\mathbb Z_9$, and
$\mathbb Z_3\times\mathbb Z_3$.
\end{Conjecture}

Horton's conjecture remained open for nearly four decades. For cyclic groups, the corresponding existence problem is whether every cyclic group $\Z_v$ of odd order $v\ge 5$, except for $\Z_5$ and $\Z_9$, admits a strong starter. Both Horton's conjecture and the corresponding cyclic existence problem have continued to be stated as conjectures or open problems in numerous subsequent surveys, monographs, and papers; see, for example, Dinitz and Stinson \cite[Section 18]{DS1992}, Stinson \cite[Section 5.4]{Stinson2004}, Dinitz \cite[Remark~55.31]{handbook}, Evans \cite[Section 16.4.5]{Evans2018}, and Mella and Pasotti \cite[Conjecture 4.2]{MP}.

In this paper, we settle Horton's conjecture and obtain a complete characterization of finite abelian groups of odd order that admit strong starters.

\begin{Theorem}\label{abelian}
Let $G$ be a nontrivial finite abelian group of odd order. Then there exists a strong starter in $G$ if and only if $G$ is not isomorphic to any of $\Z_3$, $\Z_5$, $\Z_9$, and $\Z_3\times\Z_3$.
\end{Theorem}

Strong starters have applications in the construction of various designs, including Room squares, Howell designs, Kirkman triple systems, and Kirkman squares (cf. \cite[Section 1]{ksv}). As applications of Theorem \ref{abelian}, we determine in Section \ref{sec-application}  which cycles of odd order or twice an odd order admit optimal extended irregular dominating sets. We also establish the existence of three pairwise orthogonal starters in finite abelian groups of odd order, and the existence of cyclic Room squares.

The remainder of this paper is organized as follows. Section \ref{sec-main} concerns strong starters in cyclic groups, Section \ref{sec-noncyclic} deals with strong starters in noncyclic abelian groups, and Section \ref{sec-application} presents applications of our results.

\section{Strong starters in cyclic groups}\label{sec-main}

In this section, we establish the existence of strong starters in cyclic groups. We always identify the elements in $\Z_v$ with their least nonnegative integer representatives $0,1,\ldots,v-1$ under the standard integer order. For positive integers $a,b,c$ with $a\le b$ and $a\equiv b\pmod c$, we set $[a,b]_c:=\{a+cj:0\le j\le(b-a)/c\}$, and write $[a,b]$ for $[a,b]_1$. 

To construct a strong starter in $\Z_v$ with $v\equiv3\pmod6$, we first form 12 families of unordered pairs, where the two entries in each pair depend linearly on an integer $i$, with coefficients $1$ and $-2$, respectively. We then add a few further pairs. The constants in these families and the additional pairs are chosen so that all entries are distinct and nonzero, the differences of
the pairs are pairwise distinct, and the pair sums are distinct and nonzero. A similar construction strategy was used in
\cite{LFWZ26} to construct perfect difference families.

\begin{Lemma}\label{ss}
For every integer $v\equiv3\pmod6$ with $v\ge201$, there exists a strong starter in $\Z_v$.
\end{Lemma}

\begin{proof}
Write $v=24t+r$, where $r\in\{3,9,15,21\}$. Since $v\ge201$, we have $t\ge9$ when $r=3$, and $t\ge8$ otherwise. We construct $(v-1)/2=12t+(r-1)/2$ unordered pairs in two parts. The first part consists of the $12(t-3)$ pairs:
$$S_{i,j}:=\{a_jt+i+\alpha_j,\ b_jt-2i+\beta_j\},\quad 2\le i\le t-2,\quad 1\le j\le12,$$
where
\begin{align*}
(a_1,a_2,\ldots,a_{12})&=(2,3,4,5,6,15,16,19,20,21,22,23),\\
(b_1,b_2,\ldots,b_{12})&=(2,15,19,11,9,9,13,13,11,15,19,2),
\end{align*}
and $\alpha_j,\beta_j$ are given in Table \ref{tab:parameters}. 

{\setlength{\tabcolsep}{2.5pt}
\begin{longtable}{c|cc|cc|cc|cc|cc|cc|cc|cc|cc|cc|cc|cc}\hline
$r$ & $\alpha_1$ & $\beta_1$& $\alpha_2$ & $\beta_2$& $\alpha_3$ & $\beta_3$& $\alpha_4$ & $\beta_4$& $\alpha_5$ & $\beta_5$& $\alpha_6$ & $\beta_6$& $\alpha_7$ & $\beta_7$& $\alpha_8$ & $\beta_8$& $\alpha_9$ & $\beta_9$& $\alpha_{10}$ & $\beta_{10}$& $\alpha_{11}$ & $\beta_{11}$& $\alpha_{12}$ & $\beta_{12}$ \\
\hline
$3$ & $-3$ & $1$ & $1$ & $0$ & $-1$ & $2$ & $-3$ & $3$ & $-3$ & $4$ & $-1$ & $-1$ & $5$ & $-2$ & $5$ & $3$ & $5$ & $4$ & $4$ & $3$ & $3$ & $1$ & $4$ & $-2$\\
$9$ & $1$ & $-2$ & $-1$ & $3$ & $0$ & $9$ & $2$ & $4$ & $5$ & $0$ & $8$ & $1$ & $6$ & $6$ & $6$ & $3$ & $6$ & $1$ & $10$ & $2$ & $8$ & $4$ & $8$ & $1$\\
$15$ & $0$ & $3$ & $1$ & $8$ & $1$ & $14$ & $2$ & $6$ & $2$ & $4$ & $6$ & $7$ & $9$ & $7$ & $10$ & $6$ & $11$ & $5$ & $14$ & $11$ & $11$ & $11$ & $12$ & $4$\\
$21$ & $-2$ & $2$ & $-1$ & $14$ & $6$ & $14$ & $5$ & $12$ & $2$ & $8$ & $13$ & $7$ & $13$ & $13$ & $14$ & $12$ & $18$ & $7$ & $16$ & $15$ & $19$ & $17$ & $17$ & $3$ \\\hline
\noalign{\vskip 8pt}
\caption{Values of $\alpha_j$ and $\beta_j$}\label{tab:parameters}
\end{longtable}}

The remaining $36+(r-1)/2$ pairs are listed in Table \ref{tab:remaining-pairs}. Let $\mathcal F=\{S_{i,j}:2\le i\le t-2,\,1\le j\le12\}$ and let $\mathcal F'$ be the collection of the remaining $36+(r-1)/2$ pairs. We claim that $\mathcal F\cup \mathcal F'$ is a strong starter in $\Z_v$.
 
\begin{longtable}{llll}\hline
$r=3$ & &&\\
\hline
$\{2t-4, 23t+5\}$, & $\{7t-4, 7t+6\}$, & $\{6t-4, 9t+5\}$,& $\{13t-4, 16t+6\}$, \\
$\{3, 3t-4\}$, & $\{2t-2, 17t+4\}$,& $\{11t+5, 20t+4\}$, & $\{20t+6, 23t+3\}$, \\
$\{11t+3, 20t+5\}$, &$\{11t+1, 23t+4\}$, & $\{5t-2, 23t+2\}$, & $\{7t+4, 22t+3\}$, \\
$\{1, 9t-1\}$, & $\{9t+1, 21t+5\}$, & $\{9t-3, 9t+4\}$, & $\{13t-2, 16t+5\}$, \\
$\{6t-3, 15t\}$, & $\{3t-3, 3t+1\}$, & $\{3t-2, 3t-1\}$, & $\{3t, 9t+3\}$, \\
$\{6t-2, 9t+2\}$, &$\{3t+2, 15t-2\}$, & $\{9t+6, 21t+4\}$, & $\{4t, 13t\}$, \\
$\{7t+2, 22t+4\}$, & $\{16t+4, 19t\}$, & $\{7t, 13t+5\}$, & $\{16t-2, 19t+4\}$, \\
$\{16t+1, 19t+6\}$, & $\{7t-2, 19t-1\}$, & $\{7t+1,13t+3\}$, & $\{7t-1, 19t+5\}$, \\
$\{7t-3, 13t+1\}$, &$\{13t+2, 16t+2\}$, & $\{16t, 19t+3\}$, & $\{16t-1, 19t+1\}$, \\
$\{16t+3, 19t+2\}$.&&&\\
\hline
$r=9$ & &&\\
\hline
$\{6t+6, 9t-2\}$, & $\{15t, 21t+11\}$, & $\{5t-1, 17t+11\}$, & $\{6t+5, 9t\}$, \\
$\{11t+2, 17t+9\}$, & $\{9t-1, 15t+9\}$, & $\{15t+2, 21t+10\}$, & $\{15t+8, 21t+5\}$, \\
$\{2t-4, 5t+3\}$, & $\{19t+2, 19t+7\}$, & $\{11t-1, 11t+5\}$, & $\{3, 21t+7\}$, \\
$\{4t-1, 22t+9\}$,& $\{4t, 16t+7\}$, & $\{11t+8, 23t+8\}$, & $\{4t-2, 13t+5\}$, \\
$\{4t+1, 13t+1\}$, & $\{13t+3, 19t+6\}$, & $\{13t+4, 19t+4\}$, & $\{9t+6, 24t+7\}$, \\
$\{6t+4, 9t+2\}$, & $\{11t+1, 17t+7\}$, & $\{2t+1, 20t+5\}$, & $\{21t+9, 24t+8\}$, \\
$\{9t+1, 9t+3\}$, & $\{3t, 15t+3\}$, & $\{15t+1, 15t+4\}$, & $\{11t+6, 23t+7\}$, \\
$\{2t+2, 20t+7\}$, & $\{6t+3, 9t+4\}$, & $\{1, 21t+6\}$, & $\{6t+1, 15t+5\}$, \\
$\{6t+2, 15t+7\}$, & $\{2t-2, 5t\}$, & $\{15t+6, 21t+8\}$, & $\{11t+4, 17t+5\}$, \\
$\{2t-1, 5t+2\}$, & $\{5t+1, 17t+6\}$, & $\{11t+3, 20t+6\}$, & $\{2t, 23t+9\}$.\\	
\hline
$r=15$ && &\\
\hline
$\{6, 24t+11\}$, & $\{7t+9, 16t+5\}$, & $\{5t, 17t+16\}$, & $\{5, 3t-1\}$, \\
$\{3, 9t+2\}$, & $\{4, 9t+4\}$, & $\{13t+4, 13t+11\}$, & $\{17t+13, 23t+10\}$, \\
$\{6t+1, 21t+15\}$, & $\{2, 21t+10\}$, & $\{11t+8, 23t+11\}$, & $\{5t+1, 17t+14\}$, \\
$\{2t+1, 11t+3\}$, & $\{7t+1, 7t+6\}$, & $\{11t+9, 17t+8\}$, & $\{13t+9, 13t+13\}$, \\
$\{17t+12, 20t+9\}$, & $\{7t+3, 13t+10\}$, & $\{9t+8, 24t+14\}$, & $\{9t+7, 24t+12\}$, \\
$\{16t+10, 19t+9\}$, & $\{4t, 7t+5\}$, & $\{4t+2, 16t+6\}$, & $\{9t+6, 21t+11\}$, \\
$\{7t+7, 13t+7\}$, & $\{3t+2, 21t+14\}$, & $\{6t+3, 15t+6\}$, & $\{16t+7, 19t+11\}$, \\
$\{7t+2, 13t+8\}$, & $\{7t+4, 13t+5\}$, & $\{4t+1, 13t+6\}$, & $\{16t+8, 16t+9\}$, \\
$\{3t+1, 21t+12\}$, & $\{17t+10, 23t+12\}$, & $\{1, 21t+13\}$, & $\{3t, 6t+2\}$, \\
$\{9t+5, 24t+13\}$, & $\{11t+6, 20t+10\}$, & $\{11t+4, 23t+13\}$, & $\{5t+2, 17t+9\}$, \\
$\{5t+3, 20t+12\}$, & $\{11t+5, 11t+7\}$, & $\{17t+11, 20t+11\}$.&\\
\hline
$r=21$ & &&\\
\hline
$\{11t+5, 11t+15\}$, & $\{4, 3t-3\}$, & $\{13t+17, 22t+15\}$, & $\{13t+15, 19t+12\}$, \\
$\{11t+14, 20t+13\}$, & $\{5, 9t+6\}$, & $\{11t+13, 17t+13\}$, & $\{9t+12, 24t+17\}$, \\
$\{4t-2, 13t+16\}$, & $\{17t+19, 23t+18\}$, & $\{9t+5, 9t+14\}$, & $\{4t, 7t+10\}$, \\
$\{9t+10, 24t+16\}$, & $\{13t+14, 22t+16\}$, & $\{4t+7, 19t+14\}$, & $\{11t+11, 20t+14\}$, \\
$\{7t+9, 22t+17\}$, & $\{11t+12, 17t+15\}$, & $\{4t-1, 4t+6\}$, & $\{7t+2, 7t+8\}$, \\
$\{2, 3t-2\}$, & $\{3t-1, 3t\}$, & $\{17t+12, 20t+19\}$, & $\{17t+17, 20t+16\}$, \\
$\{11t+7, 11t+10\}$, & $\{7t+1, 7t+5\}$, & $\{7t+3, 13t+12\}$, & $\{13t+13, 19t+15\}$, \\
$\{4t+5, 16t+13\}$, & $\{4t+4, 22t+20\}$, & $\{17t+14, 20t+18\}$, & $\{5t+5, 17t+16\}$, \\
$\{4t+3, 22t+18\}$, & $\{4t+2, 7t+4\}$, & $\{4t+1, 7t+6\}$, & $\{7t+7, 16t+12\}$, \\
$\{13t+10, 22t+19\}$, & $\{13t+11, 16t+14\}$, & $\{15t+13, 21t+17\}$, & $\{9t+9, 24t+19\}$, \\
$\{1, 9t+8\}$, & $\{5t+6, 20t+17\}$, & $\{11t+9, 20t+15\}$, & $\{15t+14, 24t+18\}$, \\
$\{3, 15t+12\}$, & $\{9t+7, 24t+20\}$.&&\\
\hline\noalign{\vskip 8pt}
\caption{The remaining $36+(r-1)/2$ pairs}
\label{tab:remaining-pairs}
\end{longtable}

We give the verification for the case $r=3$; the other three cases $r=9,15,21$ can be verified similarly. We first show that the pairs in $\mathcal F\cup\mathcal F'$ partition $\Z_v\setminus\{0\}$. For each fixed $j$, as $i$ ranges from $2$ to $t-2$, the two entries of $S_{i,j}=\{a_jt+i+\alpha_j,\ b_jt-2i+\beta_j\}$ take precisely the values listed in the second and third columns of Table \ref{tab:points}, respectively. The sets of values listed in these two columns are pairwise disjoint and contain only nonzero elements of $\Z_v$. Thus, the pairs in $\mathcal F$ contain $24(t-3)=v-75$ distinct nonzero elements of $\Z_v$. The remaining $37$ pairs in $\mathcal F'$ contain precisely the other $74$ nonzero elements, each occurring exactly once. 

\begin{table}[!ht]
\centering
\begin{tabular}{c|c|c}
\hline
$j$ & $a_jt+i+\alpha_j$ & $b_jt-2i+\beta_j$ \\
\hline
$1$ & $2t+i-3\in [2t-1, 3t-5]$ & $2t-2i+1\in [5, 2t-3]_2$ \\
\hline
$2$ & $3t+i+1\in [3t+3, 4t-1]$ & $15t-2i\in [13t+4, 15t-4]_2$ \\
\hline
$3$ & $4t+i-1\in [4t+1, 5t-3]$ & $19t-2i+2\in [17t+6, 19t-2]_2$ \\
\hline
$4$ & $5t+i-3\in [5t-1, 6t-5]$ & $11t-2i+3\in [9t+7, 11t-1]_2$ \\
\hline
$5$ & $6t+i-3\in [6t-1, 7t-5]$ & $9t-2i+4\in [7t+8, 9t]_2$ \\
\hline
$6$ & $15t+i-1\in [15t+1, 16t-3]$ & $9t-2i-1\in [7t+3, 9t-5]_2$ \\
\hline
$7$ & $16t+i+5\in [16t+7, 17t+3]$ & $13t-2i-2\in [11t+2, 13t-6]_2$ \\
\hline
$8$ & $19t+i+5\in [19t+7, 20t+3]$ & $13t-2i+3\in [11t+7, 13t-1]_2$ \\
\hline
$9$ & $20t+i+5\in [20t+7, 21t+3]$ & $11t-2i+4\in [9t+8, 11t]_2$ \\
\hline
$10$ & $21t+i+4\in [21t+6, 22t+2]$ & $15t-2i+3\in [13t+7, 15t-1]_2$ \\
\hline
$11$ & $22t+i+3\in [22t+5, 23t+1]$ & $19t-2i+1\in [17t+5, 19t-3]_2$ \\
\hline
$12$ & $23t+i+4\in [23t+6, 24t+2]$ & $2t-2i-2\in [2, 2t-6]_2$ \\
\hline
\end{tabular}
\caption{Elements of the pairs in $\mathcal F$ for $v=24t+3$.}
\label{tab:points}
\end{table}

We next verify the difference condition. For distinct $x,y\in\{0,1,\ldots,v-1\}$, define the \emph{positive difference} of the pair $\{x,y\}$ by $d_v(\{x,y\})=\min\{|x-y|,\,v-|x-y|\}$. Since $v$ is odd, exactly one of the two differences $x-y$ and $y-x$, viewed modulo $v$, has its least nonnegative representative in $[1,(v-1)/2]$. This representative is precisely $d_v(\{x,y\})$. Consequently, the difference condition is equivalent to requiring that the positive differences of the pairs in $\mathcal F\cup\mathcal F'$ cover $[1,(v-1)/2]$ exactly once. The positive differences of the pairs in $\mathcal F$ are listed in the second column of Table \ref{tab:d-s}, while those of the pairs in $\mathcal F'$ are listed in the second and fifth columns of Table \ref{tab:37}. All these values are pairwise distinct and belong to $[1,12t+1]$. Since their total number is $12(t-3)+37=12t+1=(v-1)/2$, they comprise precisely all the integers in $[1,(v-1)/2]$. Thus, the difference condition holds.

\begin{table}[!ht]
\centering
\begin{tabular}{c|c|c}
\hline
$j$ & Positive differences of $S_{i,j}$ & Sums of $S_{i,j}$ \\
\hline
$1$ & $3i-4\in[2,3t-10]_3$ & $4t-i-2\in[3t,4t-4]$ \\
\hline
$2$ & $12t-3i-1\in[9t+5,12t-7]_3$ & $18t-i+1\in[17t+3,18t-1]$ \\
\hline
$3$ & $9t+3i\in[9t+6,12t-6]_3$ & $23t-i+1\in[22t+3,23t-1]$ \\
\hline
$4$ & $6t-3i+6\in[3t+12,6t]_3$ & $16t-i\in[15t+2,16t-2]$ \\
\hline
$5$ & $3t-3i+7\in[13,3t+1]_3$ & $15t-i+1\in[14t+3,15t-1]$ \\
\hline
$6$ & $6t+3i\in[6t+6,9t-6]_3$ & $24t-i-2\in[23t,24t-4]$ \\
\hline
$7$ & $3t+3i+7\in[3t+13,6t+1]_3$ & $5t-i\in[4t+2,5t-2]$ \\
\hline
$8$ & $6t+3i+2\in[6t+8,9t-4]_3$ & $8t-i+5\in[7t+7,8t+3]$ \\
\hline
$9$ & $9t+3i+1\in[9t+7,12t-5]_3$ & $7t-i+6\in[6t+8,7t+4]$ \\
\hline
$10$ & $6t+3i+1\in[6t+7,9t-5]_3$ & $12t-i+4\in[11t+6,12t+2]$ \\
\hline
$11$ & $3t+3i+2\in[3t+8,6t-4]_3$ & $17t-i+1\in[16t+3,17t-1]$ \\
\hline
$12$ & $3t-3i-3\in[3,3t-9]_3$ & $t-i-1\in[1,t-3]$ \\
\hline
\end{tabular}
\caption{Positive differences and sums of the pairs in $\mathcal F$ for $v=24t+3$.}
\label{tab:d-s}
\end{table}

It remains to verify the strong condition. The sums of the pairs in $\mathcal F$ are listed in the third column of Table~\ref{tab:d-s}, and the sums of the pairs in $\mathcal F'$ are listed in the third and sixth columns of Table~\ref{tab:37}. For $t\ge9$, the intervals in the third column of Table~\ref{tab:d-s} are pairwise disjoint. Moreover, the $37$ sums listed in the third and sixth columns of Table \ref{tab:37} are pairwise distinct and lie outside the union of these intervals. Since all the sums are nonzero, the strong condition follows.

This completes the verification that $\mathcal F\cup\mathcal F'$ is a strong starter in $\Z_v$ for $r=3$. The other three cases $r=9,15,21$ can be verified similarly. For the convenience of the reader, the Python code used to verify our results is available at \cite{StrongStarterCode}.
\end{proof}

\begin{table}[!ht]
\centering
\begin{tabular}{c|c|c||c|c|c}
\hline
Pair & \shortstack{Positive\\difference} & Sum
& Pair & \shortstack{Positive\\difference} & Sum\\
\hline
$\{2t-4,23t+5\}$ & $3t-6$ & $t-2$ & $\{7t-4,7t+6\}$ & $10$ & $14t+2$ \\
\hline
$\{6t-4,9t+5\}$ & $3t+9$ & $15t+1$ & $\{13t-4,16t+6\}$ & $3t+10$ & $5t-1$ \\
\hline
$\{3,3t-4\}$ & $3t-7$ & $3t-1$ & $\{2t-2,17t+4\}$ & $9t-3$ & $19t+2$ \\
\hline
$\{11t+5,20t+4\}$ & $9t-1$ & $7t+6$ & $\{20t+6,23t+3\}$ & $3t-3$ & $19t+6$ \\
\hline
$\{11t+3,20t+5\}$ & $9t+2$ & $7t+5$ & $\{11t+1,23t+4\}$ & $12t$ & $10t+2$ \\
\hline
$\{5t-2,23t+2\}$ & $6t-1$ & $4t-3$ & $\{7t+4,22t+3\}$ & $9t+4$ & $5t+4$ \\
\hline
$\{1,9t-1\}$ & $9t-2$ & $9t$ & $\{9t+1,21t+5\}$ & $12t-1$ & $6t+3$ \\
\hline
$\{9t-3,9t+4\}$ & $7$ & $18t+1$ & $\{13t-2,16t+5\}$ & $3t+7$ & $5t$ \\
\hline
$\{6t-3,15t\}$ & $9t+3$ & $21t-3$ & $\{3t-3,3t+1\}$ & $4$ & $6t-2$ \\
\hline
$\{3t-2,3t-1\}$ & $1$ & $6t-3$ & $\{3t,9t+3\}$ & $6t+3$ & $12t+3$ \\
\hline
$\{6t-2,9t+2\}$ & $3t+4$ & $15t$ & $\{3t+2,15t-2\}$ & $12t-4$ & $18t$ \\
\hline
$\{9t+6,21t+4\}$ & $12t-2$ & $6t+7$ & $\{4t,13t\}$ & $9t$ & $17t$ \\
\hline
$\{7t+2,22t+4\}$ & $9t+1$ & $5t+3$ & $\{16t+4,19t\}$ & $3t-4$ & $11t+1$ \\
\hline
$\{7t,13t+5\}$ & $6t+5$ & $20t+5$ & $\{16t-2,19t+4\}$ & $3t+6$ & $11t-1$ \\
\hline
$\{16t+1,19t+6\}$ & $3t+5$ & $11t+4$ & $\{7t-2,19t-1\}$ & $12t+1$ & $2t-6$ \\
\hline
$\{7t+1,13t+3\}$ & $6t+2$ & $20t+4$ & $\{7t-1,19t+5\}$ & $12t-3$ & $2t+1$ \\
\hline
$\{7t-3,13t+1\}$ & $6t+4$ & $20t-2$ & $\{13t+2,16t+2\}$ & $3t$ & $5t+1$ \\
\hline
$\{16t,19t+3\}$ & $3t+3$ & $11t$ & $\{16t-1,19t+1\}$ & $3t+2$ & $11t-3$ \\
\hline
$\{16t+3,19t+2\}$ & $3t-1$ & $11t+2$ &&&\\
\hline
\end{tabular}
\caption{Positive differences and sums of the $37$ pairs in $\mathcal F'$ for $v=24t+3$.}
\label{tab:37}
\end{table}

Combining Lemma \ref{ss} with Lemmas~\ref{999} and~\ref{known}, we determine the existence spectrum of strong starters in cyclic groups.

\begin{Proposition}\label{main}
For any odd integer $v\ge 7$ and $v\neq 9$, there exists a strong starter in $\Z_v$.
\end{Proposition}

\section{Strong starters in abelian groups}\label{sec-noncyclic}

In this section, we establish the existence of strong starters in abelian groups. Horton \cite[Theorem 2.14]{H} showed that Conjecture \ref{conj} can be reduced to the existence of strong starters in the following five families of groups. We write $\Z_3^2=\Z_3\times\Z_3$ for the direct product of two copies of the additive group $\Z_3$.

\begin{Lemma}{\rm \cite{H}}\label{lem:Hor-abe}
Conjecture \ref{conj} holds if and only if strong starters exist in each of the following families of abelian groups:
\begin{enumerate}
\item[$(1)$] $\Z_{3p}$, where $p>5$ is a prime;
\item[$(2)$] $\Z_{9p}$, where $p>5$ is a prime;
\item[$(3)$] $\Z_3^2\times\Z_p$, where $p>5$ is a prime;
\item[$(4)$] $\Z_{3^a}$, where $a\ge 7$;
\item[$(5)$] $\Z_{5\cdot 3^a}$, where $a\ge 4$.
\end{enumerate}
\end{Lemma}

By Proposition \ref{main}, strong starters exist in all the cyclic groups listed in Lemma \ref{lem:Hor-abe}. Thus, it remains to
consider the noncyclic family $\Z_3^2\times\Z_p$, where $p>5$ is a prime.

A useful tool for constructing strong starters is the notion of a frame starter. Let $(G,+)$ be an abelian group of order $g$, and let $H$ be a subgroup of $G$ of order $h$, where $g-h$ is even. A \emph{frame starter} in $G\setminus H$ is a set of unordered pairs $\mathcal S=\{\{x_i,y_i\}:1\le i\le(g-h)/2\}$ satisfying the following conditions:
\begin{enumerate}
\item[(i)] $\{x_i:1\le i\le(g-h)/2\}\cup \{y_i:1\le i\le(g-h)/2\} = G\setminus H$;
\item[(ii)] $\{\pm(x_i-y_i):1\le i\le(g-h)/2\} = G\setminus H$.
\end{enumerate}
A frame starter is \emph{strong} if, in addition, the sums $x_i+y_i$, for $1\le i\le(g-h)/2$, are pairwise distinct and belong to $G\setminus H$. When $H=\{0_G\}$, a frame starter or strong frame starter in $G\setminus H$ is just a starter or strong starter in $G$, respectively.

The notion of a partial starter appeared in \cite{AG78}, and frame starters were subsequently introduced in \cite{DS80Ars}.
For further information on frame starters, we refer the reader to~\cite{Stinson2024}. For our purpose, we need the following elementary result, whose proof is straightforward and is omitted.

\begin{Lemma}\label{lem-frame-fill}
Let $H$ be a subgroup of an abelian group $G$ of odd order. If there is a strong frame starter in $G\setminus H$ and a strong starter in $H$, then there is a strong starter in $G$.
\end{Lemma}

By Lemma \ref{lem-frame-fill}, to construct a strong starter in $\Z_3^2\times\Z_p$, it suffices to find a subgroup $H$ of
$\Z_3^2\times\Z_p$ that admits a strong starter and for which there is a strong frame starter in $(\Z_3^2\times\Z_p)\setminus H$. We take $H=\{(0,0)\}\times\Z_p$, which is isomorphic to $\Z_p$. By Proposition \ref{main}, $H$ admits a strong starter. Thus, it remains to construct a strong frame starter in $(\Z_3^2\times\Z_p)\setminus(\{(0,0)\}\times\Z_p)$.

Throughout the remainder of this section, for a prime $p$, we denote by $\Z_p^*$ the multiplicative group of nonzero elements of $\Z_p$. An element of $\Z_3^2\times\Z_p$ is written as $(x,t)$, or as $(a,b,t)$ when $x=(a,b)\in\Z_3^2$. Set
$$
D=\{(1,0),(0,1),(1,1),(2,1)\}.
$$
Then $D\cup(-D)=\Z_3^2\setminus\{(0,0)\}$.

\begin{Lemma}\label{lem-frame-even}
Let $p>5$ be a prime. If $2$ has even multiplicative order in $\Z_p^*$, then there is a strong frame starter in $(\Z_3^2\times\Z_p)\setminus(\{(0,0)\}\times\Z_p)$.
\end{Lemma}

\begin{proof}
Since $2$ has even multiplicative order, the cyclic subgroup $\langle 2\rangle$ of $\Z_p^*$ has even order. For each coset of
$\langle 2\rangle$ in $\Z_p^*$, choose a representative $t_0$, taking $t_0=1$ for the coset $\langle 2\rangle$, and define $f$
on this coset by
$$
f(2^i t_0)=(-1)^i,\qquad i=0,1,\ldots,|\langle 2\rangle|-1.
$$
Since $|\langle 2\rangle|$ is even, this defines a function $f:\Z_p^*\to\{1,-1\}$ satisfying
$$
f(2t)=-f(t),\qquad f(1)=1.
$$
For each $d\in D$, choose $x_d,y_d\in\Z_3^2$ as follows:
$$
\begin{array}{ccc|ccc}\hline
\multicolumn{3}{c|}{f(3^{-1})=1}
& \multicolumn{3}{c}{f(3^{-1})=-1}\\
\hline
d & x_d & y_d & d & x_d & y_d\\
\hline
(1,0) & (0,1) & (2,2) & (1,0) & (0,1) & (0,2)\\
(0,1) & (1,1) & (2,1) & (0,1) & (1,1) & (2,1)\\
(1,1) & (1,2) & (1,0) & (1,1) & (2,0) & (1,0)\\
(2,1) & (1,0) & (0,2) & (2,1) & (0,2) & (2,2)\\
\hline
\end{array}
$$
Let
\begin{align*}
\mathcal A={}&\bigl\{\{(x_d,t),(x_d+d,2t)\}:d\in D,\ t\in\Z_p^*\setminus\{1\},\ f(t)=1\bigr\}\\
&\quad\cup\bigl\{\{(y_d,t),(y_d+d,2t)\}:d\in D,\ t\in\Z_p^*\setminus\{1\},\ f(t)=-1\bigr\}.
\end{align*}
Let $\mathcal B$ consist of the 8 pairs listed in Table \ref{tab-frame-even-replacement} according to the value of $f(3^{-1})$. We claim that $\mathcal A\cup \mathcal B$ is a strong frame starter in $(\Z_3^2\times\Z_p)\setminus(\{(0,0)\}\times\Z_p)$.

\begin{table}[H]
\centering\setstretch{1.06}
\begin{tabular}{c|c|c|c|c|c}\hline
\multicolumn{3}{c|}{$f(3^{-1})=1$} & \multicolumn{3}{c}{$f(3^{-1})=-1$}\\\hline
Pair & Differences & Sum & Pair & Differences & Sum\\\hline
$\{(2,1,0),(0,1,1)\}$ & $\pm(1,0,1)$ & $(2,2,1)$ & $\{(0,2,1),(2,0,2)\}$ & $\pm(2,1,1)$ & $(2,2,3)$\\
$\{(1,0,1),(1,1,2)\}$ & $\pm(0,1,1)$ & $(2,1,3)$ & $\{(2,0,1),(0,1,2)\}$ & $\pm(1,1,1)$ & $(2,1,3)$\\
$\{(1,2,2),(0,1,2)\}$ & $\pm(1,1,0)$ & $(1,0,4)$ & $\{(2,1,0),(0,1,1)\}$ & $\pm(1,0,1)$ & $(2,2,1)$\\
$\{(1,2,1),(2,0,2)\}$ & $\pm(1,1,1)$ & $(0,2,3)$ & $\{(1,1,2),(1,2,2)\}$ & $\pm(0,1,0)$ & $(2,0,4)$\\
$\{(2,0,0),(1,1,1)\}$ & $\pm(2,1,1)$ & $(0,1,1)$ & $\{(1,0,0),(1,1,1)\}$ & $\pm(0,1,1)$ & $(2,1,1)$\\
$\{(1,0,0),(0,1,0)\}$ & $\pm(2,1,0)$ & $(1,1,0)$ & $\{(2,0,0),(1,1,0)\}$ & $\pm(2,1,0)$ & $(0,1,0)$\\
$\{(1,1,0),(1,2,0)\}$ & $\pm(0,1,0)$ & $(2,0,0)$ & $\{(0,1,0),(1,2,0)\}$ & $\pm(1,1,0)$ & $(1,0,0)$\\
$\{(0,2,0),(2,2,0)\}$ & $\pm(1,0,0)$ & $(2,1,0)$ & $\{(0,2,0),(2,2,0)\}$ & $\pm(1,0,0)$ & $(2,1,0)$\\
\hline\end{tabular}
\caption{Differences and sums of the pairs in $\mathcal B$ in Lemma \ref{lem-frame-even}}
\label{tab-frame-even-replacement}
\end{table}

We first show that the pairs in $\mathcal A\cup\mathcal B$ partition $(\Z_3^2\setminus\{(0,0)\})\times\Z_p$. Let $r\in\Z_p^*\setminus\{1,2\}$. The elements of $\mathcal A$ with last coordinate $r$ are
$$
\{(x_d,r):d\in D,\ f(r)=1\}\cup\{(y_d,r):d\in D,\ f(r)=-1\}
$$
together with
$$
\{(x_d+d,r):d\in D,\ f(2^{-1}r)=1\}\cup\{(y_d+d,r):d\in D,\ f(2^{-1}r)=-1\}.
$$
Since $f(2^{-1}r)=-f(r)$, the above elements are precisely
$$
\begin{cases}
\{(x_d,r):d\in D\}\cup\{(y_d+d,r):d\in D\},
&\text{if } f(r)=1;\\
\{(y_d,r):d\in D\}\cup\{(x_d+d,r):d\in D\},
&\text{if } f(r)=-1.
\end{cases}
$$
Since $\{x_d:d\in D\}\cup\{y_d+d:d\in D\}=\{y_d:d\in D\}\cup\{x_d+d:d\in D\}=\Z_3^2\setminus\{(0,0)\}$, 
the elements of $\mathcal A$ with last coordinate $r$ form precisely $(\mathbb Z_3^2\setminus\{(0,0)\})\times\{r\}$. It remains to examine $r=0,1,2$. The eight elements in $\mathcal B$ with last coordinate $0$ cover $(\mathbb Z_3^2\setminus\{(0,0)\})\times\{0\}$. Since $f(1)=1$ and $f(2)=-1$, the elements of $\mathcal A$ with last coordinate $1$ and $2$ are, respectively,
$$
\{(y_d+d,1):d\in D\}
\quad\text{and}\quad
\{(y_d,2):d\in D\}.
$$
By Table~\ref{tab-frame-even-replacement}, the $\Z_3^2$-components of the elements of $\mathcal B$ with last coordinate $1$ and $2$ are, respectively,
$$
\begin{cases}
\{(0,1),(1,0),(1,2),(1,1)\}
\text{ and }\{(1,1),(1,2),(0,1),(2,0)\},
&\text{if }f(3^{-1})=1,\\
\{(0,2),(2,0),(0,1),(1,1)\}
\text{ and }\{(2,0),(0,1),(1,1),(1,2)\},
&\text{if }f(3^{-1})=-1.
\end{cases}
$$
In either case, the elements of $\mathcal B$ complement those of $\mathcal A$ in the two corresponding sets $(\mathbb Z_3^2\setminus\{(0,0)\})\times\{1,2\}$. Thus $\mathcal A\cup\mathcal B$ partitions the required set.

We next verify the difference condition. The differences of the pairs in $\mathcal A$ are of the form $(d,t)$ or $(-d,-t)$, where
$d\in D$ and $t\in\Z_p^*\setminus\{1\}$. Hence, they are precisely the elements of
$
\bigl(D\times(\Z_p^*\setminus\{1\})\bigr)
\cup
\bigl((-D)\times(\Z_p^*\setminus\{-1\})\bigr).
$
By Table~\ref{tab-frame-even-replacement}, the differences of the pairs in $\mathcal B$ are precisely the elements of
$(D\times\{0,1\})\cup((-D)\times\{0,-1\})$. Since $D\cup(-D)=\Z_3^2\setminus\{(0,0)\}$, the differences of the pairs in $\mathcal A\cup\mathcal B$ comprise every element of $(\Z_3^2\setminus\{(0,0)\})\times\Z_p$ exactly once.

It remains to verify the strong condition. For each $t\in\Z_p^*\setminus\{1\}$, the sums of the pairs in
$\mathcal A$ indexed by $t$ are
$$
\{(2x_d+d,3t):d\in D\}\quad\text{if }f(t)=1,
$$
and
$$
\{(2y_d+d,3t):d\in D\}\quad\text{if }f(t)=-1.
$$
Note that $\{2x_d+d:d\in D\}\cup\{2y_d+d:d\in D\}=\Z_3^2\setminus\{(0,0)\}$. Thus, for each fixed $t\in\Z_p^*\setminus\{1\}$, the sums of the pairs in $\mathcal A$ indexed by $t$ are distinct, and none of these sums has $(0,0)$ as its $\Z_3^2$-component. Since $3\ne0$ in $\Z_p$, the sums corresponding to distinct values of $t$ are also distinct. Hence all sums of the pairs in $\mathcal A$ are pairwise distinct.  By Table~\ref{tab-frame-even-replacement}, the sums of the pairs in $\mathcal B$ are pairwise distinct, and none of these sums has $(0,0)$ as its $\Z_3^2$-component.

Finally, we show that no sum of a pair in $\mathcal A$ coincides with a sum of a pair in $\mathcal B$. The last coordinates of the
sums in $\mathcal B$ belong to $\{0,1,3,4\}$, whereas those of the sums in $\mathcal A$ are of the form $3t$, where $t\in\Z_p^*\setminus\{1\}$. Since $p>5$, the elements $0,1,3,4$ are distinct in $\Z_p$. Thus a coincidence is possible only for $t=3^{-1}$ or $t=4\cdot3^{-1}$. Since $f(4t)=f(t)$, the $\Z_3^2$-components of the sums in $\mathcal A$ corresponding to $t=3^{-1}$ and $t=4\cdot3^{-1}$ are the same, namely,
$$
\begin{cases}
\{(1,2),(2,0),(0,2),(1,1)\},
&\text{if }f(3^{-1})=1,\\
\{(1,1),(1,0),(0,1),(0,2)\},
&\text{if }f(3^{-1})=-1.
\end{cases}
$$
On the other hand, the $\Z_3^2$-components of the sums of the pairs in $\mathcal B$ with last coordinates $1$ and $4$ are, respectively,
$$
\begin{cases}
\{(2,2),(0,1)\}\text{ and }\{(1,0)\},
&\text{if }f(3^{-1})=1,\\
\{(2,2),(2,1)\}\text{ and }\{(2,0)\},
&\text{if }f(3^{-1})=-1.
\end{cases}
$$
These sets are disjoint in either case. Consequently, all sums of pairs in $\mathcal A\cup\mathcal B$ are pairwise distinct, and none of them has $(0,0)$ as its $\Z_3^2$-component. Thus the strong condition holds.
\end{proof}

\begin{Lemma}\label{lem-frame-odd}
Let $p>7$ be a prime. If $2$ has odd multiplicative order in $\Z_p^*$, then there is a strong frame starter in $(\Z_3^2\times\Z_p)\setminus(\{(0,0)\}\times\Z_p)$.
\end{Lemma}

\begin{proof}
Since $2$ has odd multiplicative order in $\Z_p^*$, we have $-1\notin\langle 2\rangle$. Consequently, for every $t\in\Z_p^*$, the cosets $t\langle 2\rangle$ and $-t\langle 2\rangle$ are distinct. Thus, the cosets of $\langle 2\rangle$ in $\Z_p^*$ can be partitioned into pairs of the form $\{t\langle 2\rangle,-t\langle 2\rangle\}$. Choose one coset from each pair, choosing $\langle 2\rangle$ from the pair $\{\langle 2\rangle,-\langle 2\rangle\}$. Define a function $g:\Z_p^*\to\{1,-1\}$ to take the value $1$ on each chosen coset and the value $-1$ on the other coset in the same pair. Then
$$
g(2t)=g(t),\qquad g(-t)=-g(t),\qquad g(1)=1.
$$
Define the permutation $\tau$ on $\Z_3^2$ by $\tau(a,b)=(-b,a)$. Let
$$
E=\{(1,0),(2,0),(1,1),(2,2)\},\qquad
C=\{(0,1),(2,1),(0,2),(1,2)\}.
$$
Then
$$
\Z_3^2\setminus\{(0,0)\}=E\cup C,\qquad -E=E,\qquad -C=C,\qquad \tau(E)=C,\qquad \tau(C)=E.
$$
For each $t\in\Z_p^*$, let $R_t=E$ if $g(t)=1$, and let $R_t=C$ otherwise. Thus
$$
R_{2t}=R_t,\qquad
R_{-t}=(\Z_3^2\setminus\{(0,0)\})\setminus R_t,\qquad
-R_t=R_t.
$$
Let
$$
\mathcal A=
\bigl\{\{(x,t),(\tau x,2t)\}:
t\in\Z_p^*\setminus\{1,-1\},\ x\in R_t\bigr\},
$$
and let $\mathcal B$ consist of the 12 pairs in Table~\ref{tab-frame-odd-replacement} chosen according to $R_{3^{-1}}$. We claim that $\mathcal A\cup \mathcal B$ is a strong frame starter in $(\Z_3^2\times\Z_p)\setminus(\{(0,0)\}\times\Z_p)$.

We first show that the pairs in $\mathcal A\cup\mathcal B$ partition $(\Z_3^2\setminus\{(0,0)\})\times\Z_p$. Let $r\in\Z_p^*\setminus\{1,-1,2,-2\}$. The elements with last coordinate $r$ that occur in pairs of
$\mathcal A$ are 
$
(R_r\times\{r\})\cup(\tau(R_{2^{-1}r})\times\{r\}).
$
Since $R_{2^{-1}r}=R_r$ and $\tau$ interchanges $E$ and $C$, these elements are precisely
$
(R_r\cup\tau(R_r))\times\{r\}
=(\Z_3^2\setminus\{(0,0)\})\times\{r\}.
$
For last coordinates $1,2,-1,-2$, respectively, the elements covered by pairs in $\mathcal A$ are $C\times\{1\}$, $E\times\{2\}$, $E\times\{-1\}$ and $C\times\{-2\}$. Consequently, the elements not covered by $\mathcal A$ are precisely
$$
\begin{aligned}
&((\Z_3^2\setminus\{(0,0)\})\times\{0\})
\cup(E\times\{1\})
\cup(C\times\{2\})
\cup(C\times\{-1\})
\cup(E\times\{-2\}).
\end{aligned}
$$
By Table~\ref{tab-frame-odd-replacement}, the pairs in $\mathcal B$ cover these elements exactly once. Thus $\mathcal A\cup\mathcal B$ partitions $(\Z_3^2\setminus\{(0,0)\})\times\Z_p$.

\begin{table}[H]
\centering\setstretch{1.02}
\begin{tabular}{c|c|c|c|c|c}\hline
\multicolumn{3}{c|}{$R_{3^{-1}}=E$} & \multicolumn{3}{c}{$R_{3^{-1}}=C$}\\\hline
Pair & Differences & Sum & Pair & Differences & Sum\\\hline
$\{(1,0,-2),(2,2,-2)\}$ & $\pm(2,1,0)$ & $(0,2,-4)$ & $\{(1,0,-2),(1,1,-2)\}$ & $\pm(0,1,0)$ & $(2,1,-4)$\\
$\{(2,0,-2),(2,1,-1)\}$ & $\pm(0,1,1)$ & $(1,1,-3)$ & $\{(2,0,-2),(2,1,-1)\}$ & $\pm(0,1,1)$ & $(1,1,-3)$\\
$\{(1,1,-2),(0,1,-1)\}$ & $\pm(2,0,1)$ & $(1,2,-3)$ & $\{(2,2,-2),(0,1,-1)\}$ & $\pm(1,2,1)$ & $(2,0,-3)$\\
$\{(0,2,-1),(2,0,0)\}$ & $\pm(2,1,1)$ & $(2,2,-1)$ & $\{(1,2,-1),(0,1,0)\}$ & $\pm(2,2,1)$ & $(1,0,-1)$\\
$\{(1,2,-1),(2,2,0)\}$ & $\pm(1,0,1)$ & $(0,1,-1)$ & $\{(0,2,-1),(2,2,0)\}$ & $\pm(2,0,1)$ & $(2,1,-1)$\\
$\{(1,2,0),(1,1,1)\}$ & $\pm(0,2,1)$ & $(2,0,1)$ & $\{(1,0,0),(2,1,0)\}$ & $\pm(1,1,0)$ & $(0,1,0)$\\
$\{(2,2,1),(0,1,2)\}$ & $\pm(1,2,1)$ & $(2,0,3)$ & $\{(1,1,0),(0,2,0)\}$ & $\pm(2,1,0)$ & $(1,0,0)$\\
$\{(2,0,1),(1,2,2)\}$ & $\pm(2,2,1)$ & $(0,2,3)$ & $\{(2,0,0),(1,1,1)\}$ & $\pm(2,1,1)$ & $(0,1,1)$\\
$\{(2,1,2),(0,2,2)\}$ & $\pm(1,1,0)$ & $(2,0,4)$ & $\{(1,2,0),(2,0,1)\}$ & $\pm(1,1,1)$ & $(0,2,1)$\\
$\{(1,0,0),(1,1,0)\}$ & $\pm(0,1,0)$ & $(2,1,0)$ & $\{(1,0,1),(1,2,2)\}$ & $\pm(0,2,1)$ & $(2,2,3)$\\
$\{(0,1,0),(2,1,0)\}$ & $\pm(1,0,0)$ & $(2,2,0)$ & $\{(2,2,1),(0,2,2)\}$ & $\pm(1,0,1)$ & $(2,1,3)$\\
$\{(0,2,0),(1,0,1)\}$ & $\pm(1,1,1)$ & $(1,2,1)$ & $\{(0,1,2),(2,1,2)\}$ & $\pm(1,0,0)$ & $(2,2,4)$\\
\hline\end{tabular}
\caption{Differences and sums of the pairs in $\mathcal B$ in Lemma \ref{lem-frame-odd}}
\label{tab-frame-odd-replacement}
\end{table}

We next verify the difference condition. For convenience, let $\iota$ denote the identity permutation on $\Z_3^2$. 
The maps $\tau-\iota$ and $\tau+\iota$, defined by 
$$(\tau-\iota)(a,b):=\tau(a,b)-\iota(a,b)=(-b-a,a-b)$$
and 
$$(\tau+\iota)(a,b):=\tau(a,b)+\iota(a,b)=(a-b,a+b)$$
are both permutations on $\Z_3^2$. The differences of the pairs in $\mathcal A$ are precisely the elements of
$$
\begin{aligned}
&\{(d,t):d\in(\tau-\iota)(R_t),\
t\in\Z_p^*\setminus\{1,-1\}\}
{}\cup{}
\{(-d,-t):d\in(\tau-\iota)(R_t),\
t\in\Z_p^*\setminus\{1,-1\}\}\\
=&\{(d,t):d\in(\tau-\iota)(R_t),\
t\in\Z_p^*\setminus\{1,-1\}\}
{}\cup{}
\{(d,t):d\in-(\tau-\iota)(R_{-t}),\
t\in\Z_p^*\setminus\{1,-1\}\}.
\end{aligned}
$$
Since $R_{-t}\in\{E,C\}$, and both $(\tau-\iota)(E)$ and $(\tau-\iota)(C)$ are invariant under negation, we have $-(\tau-\iota)(R_{-t})=(\tau-\iota)(R_{-t})$. Since $\tau-\iota$ is a permutation on $\Z_3^2$ and $R_t\cup R_{-t}=\Z_3^2\setminus\{(0,0)\}$, these differences of the pairs in $\mathcal A$ comprise precisely $(\Z_3^2\setminus\{(0,0)\})\times(\Z_p^*\setminus\{1,-1\})$. By Table~\ref{tab-frame-odd-replacement}, the differences of the pairs in $\mathcal B$ comprise precisely $(\Z_3^2\setminus\{(0,0)\})\times\{0,1,-1\}$. Thus, the differences of the pairs in $\mathcal A\cup\mathcal B$ contain each element of $(\Z_3^2\setminus\{(0,0)\})\times\Z_p$ exactly once.

Finally, we verify the strong condition. For each pair $\{(x,t),(\tau x,2t)\}\in\mathcal A$, its sum is $((\tau+\iota)(x),3t)$.
Since $\tau+\iota$ is a permutation on $\Z_3^2$ fixing $(0,0)$, and $x\ne(0,0)$, the sums of the pairs in $\mathcal A$ with a fixed value of $t$ are distinct, and none of these sums has $(0,0)$ as its $\Z_3^2$-component. Since $3\ne0$ in $\Z_p$, the
sums corresponding to distinct values of $t$ are also distinct. By Table~\ref{tab-frame-odd-replacement}, the sums of the pairs in $\mathcal B$ are pairwise distinct, and none of these sums has $(0,0)$ as its $\Z_3^2$-component. 

It remains to show that no sum of a pair in $\mathcal B$ coincides with a sum of a pair in $\mathcal A$. The third coordinates of the sums of the pairs in $\mathcal A$ are precisely the elements of $\Z_p^*\setminus\{3,-3\}$, whereas those of the sums of the pairs in $\mathcal B$ belong to $\{0,\pm1,\pm3,\pm4\}$. Since $p\ge 11$, the elements $0,\pm1,\pm3,\pm4$ are pairwise distinct in $\Z_p$. Thus, only the third coordinates $\pm1,\pm4$ need to be considered. By the definition of $\tau+\iota$, we have
$$
\begin{aligned}
(\tau+\iota)(E)&=\{(1,1),(2,2),(0,2),(0,1)\},\\
(\tau+\iota)(C)&=\{(2,1),(1,0),(1,2),(2,0)\}.
\end{aligned}
$$
These two sets are disjoint and their union is $\Z_3^2\setminus\{(0,0)\}$. Since $R_{2t}=R_t$ and $R_{-t}$ is the complement of $R_t$ in $\Z_3^2\setminus\{(0,0)\}$, the $\Z_3^2$-components of the sums of pairs in $\mathcal A$ with third coordinates $1$ and $4$ belong to $(\tau+\iota)(R_{3^{-1}})$, while those with third coordinates $-1$ and $-4$ belong to $(\tau+\iota)(R_{-3^{-1}})$. By Table~\ref{tab-frame-odd-replacement}, for each of these four third coordinates, the $\Z_3^2$-components of the corresponding sums in $\mathcal B$ belong to the other one of $(\tau+\iota)(E)$ and $(\tau+\iota)(C)$. Thus, the sums of pairs in $\mathcal A$ and $\mathcal B$ are mutually distinct.
\end{proof}

\begin{Lemma}\label{lem-frame-seven}
There exists a strong frame starter in $(\Z_3^2\times\Z_7)\setminus(\{(0,0)\}\times\Z_7)$.
\end{Lemma}

\begin{proof}
Let $\mathcal S$ consist of the 28 pairs in Table \ref{tab-frame-seven}. Their elements cover $(\Z_3^2\setminus\{(0,0)\})\times\Z_7$ exactly once, and the differences cover $(\Z_3^2\setminus\{(0,0)\})\times\Z_7$ exactly once. Moreover, the sums are pairwise distinct and belong to $(\Z_3^2\setminus\{(0,0)\})\times\Z_7$. Thus $\mathcal S$ is a strong frame starter.
\end{proof}

\begin{table}[H]
\centering\setstretch{1.02}
\begin{tabular}{c|c|c|c|c|c}\hline
Pair & Differences & Sum & Pair & Differences & Sum\\\hline
$\{(1,0,0),(1,1,0)\}$ & $\pm(0,1,0)$ & $(2,1,0)$ & $\{(0,1,1),(1,2,3)\}$ & $\pm(1,1,2)$ & $(1,0,4)$\\
$\{(1,0,1),(0,2,4)\}$ & $\pm(1,1,4)$ & $(1,2,5)$ & $\{(0,1,2),(2,1,0)\}$ & $\pm(1,0,2)$ & $(2,2,2)$\\
$\{(1,0,2),(1,1,6)\}$ & $\pm(0,1,4)$ & $(2,1,1)$ & $\{(0,1,3),(2,2,4)\}$ & $\pm(2,1,1)$ & $(2,0,0)$\\
$\{(1,0,3),(2,2,0)\}$ & $\pm(2,1,3)$ & $(0,2,3)$ & $\{(0,1,4),(2,1,3)\}$ & $\pm(1,0,1)$ & $(2,2,0)$\\
$\{(1,0,4),(2,2,5)\}$ & $\pm(2,1,6)$ & $(0,2,2)$ & $\{(0,1,5),(2,1,6)\}$ & $\pm(1,0,6)$ & $(2,2,4)$\\
$\{(1,0,5),(0,1,0)\}$ & $\pm(2,1,2)$ & $(1,1,5)$ & $\{(0,1,6),(1,2,4)\}$ & $\pm(1,1,5)$ & $(1,0,3)$\\
$\{(1,0,6),(1,2,0)\}$ & $\pm(0,1,6)$ & $(2,2,6)$ & $\{(1,1,1),(0,2,6)\}$ & $\pm(2,1,5)$ & $(1,0,0)$\\
$\{(2,0,0),(1,2,1)\}$ & $\pm(1,1,6)$ & $(0,2,1)$ & $\{(1,1,2),(2,1,2)\}$ & $\pm(1,0,0)$ & $(0,2,4)$\\
$\{(2,0,1),(2,2,3)\}$ & $\pm(0,1,5)$ & $(1,2,4)$ & $\{(1,1,3),(1,2,5)\}$ & $\pm(0,1,2)$ & $(2,0,1)$\\
$\{(2,0,2),(2,2,1)\}$ & $\pm(0,1,1)$ & $(1,2,3)$ & $\{(1,1,5),(2,1,1)\}$ & $\pm(1,0,3)$ & $(0,2,6)$\\
$\{(2,0,3),(1,2,2)\}$ & $\pm(1,1,1)$ & $(0,2,5)$ & $\{(2,1,4),(0,2,0)\}$ & $\pm(1,1,3)$ & $(2,0,4)$\\
$\{(2,0,4),(1,1,4)\}$ & $\pm(2,1,0)$ & $(0,1,1)$ & $\{(2,1,5),(0,2,5)\}$ & $\pm(1,1,0)$ & $(2,0,3)$\\
$\{(2,0,5),(2,2,2)\}$ & $\pm(0,1,3)$ & $(1,2,0)$ & $\{(0,2,1),(1,2,6)\}$ & $\pm(1,0,5)$ & $(1,1,0)$\\
$\{(2,0,6),(0,2,2)\}$ & $\pm(2,1,4)$ & $(2,2,1)$ & $\{(0,2,3),(2,2,6)\}$ & $\pm(1,0,4)$ & $(2,1,2)$\\
\hline\end{tabular}
\caption{A strong frame starter in $(\Z_3^2\times\Z_7)\setminus(\{(0,0)\}\times\Z_7)$.}
\label{tab-frame-seven}
\end{table}

\begin{Proposition}\label{thm-noncyclic-frame}
For every prime $p>5$, there exists a strong frame starter in $(\Z_3^2\times\Z_p)\setminus(\{(0,0)\}\times\Z_p)$.
\end{Proposition}

\begin{proof}
The case $p=7$ follows from Lemma~\ref{lem-frame-seven}, and the case $p>7$ follows from Lemmas~\ref{lem-frame-even} and~\ref{lem-frame-odd} according to the parity of the multiplicative order of $2$ modulo $p$.
\end{proof}

\begin{Corollary}\label{noncyclic}
For every prime $p>5$, there exists a strong starter in $\Z_3^2\times\Z_p$.
\end{Corollary}

\begin{proof}
Take $G=\Z_3^2\times\Z_p$ and $H=\{(0,0)\}\times\Z_p$. By Proposition \ref{thm-noncyclic-frame}, there exists a strong frame starter in $G\setminus H$. Since $p>5$, Proposition \ref{main} provides a strong starter in $H\cong\Z_p$, and the conclusion follows from Lemma~\ref{lem-frame-fill}.
\end{proof}

Now one can combine Proposition \ref{main}, Lemma \ref{lem:Hor-abe} and Corollary \ref{noncyclic} to complete the proof of Theorem \ref{abelian}.

\section{Applications}\label{sec-application}

In this section, we present two applications of strong starters in finite abelian groups of odd order. The first concerns extended irregular domination in cycles, and the second concerns orthogonal starters and cyclic Room squares.

\subsection{Extended irregular domination of cycles}

Mella and Pasotti \cite{MP} introduced the notion of an extended irregular dominating set of a graph, allowing precisely one vertex labeled with $0$ to cover itself, and proved that an odd cycle admits an optimal extended irregular dominating set if and only if the corresponding cyclic group admits a strong starter.

Let $\Gamma=(V,E)$ be a finite connected simple graph, and write $d_\Gamma(u,v)$ for the length of a shortest path between $u,v\in V$. A vertex $s$ labeled with a nonnegative integer $\ell$ \emph{covers} a vertex $u$ if $d_\Gamma(u,s)=\ell$. An \emph{extended irregular dominating set} of $\Gamma$ is a subset $S\subseteq V$ whose vertices are assigned distinct nonnegative integer labels, one of which is $0$, such that every vertex $u\in V$ is covered by some vertex $s\in S$. If such a set exists, its minimum cardinality is called the
\emph{extended irregular domination number} $\gamma_e(\Gamma)$, and a set attaining this minimum is called \emph{optimal}.

For the cycle $C_n$ of length $n$, identify its vertex set with $\Z_n$, where each vertex $i$ is adjacent to $i-1$ and $i+1$. A vertex $i$ with label $\ell$, where $1\le\ell<n/2$, covers the two vertices $i-\ell$ and $i+\ell$. If $n$ is even, a vertex with label
$n/2$ covers its antipodal vertex. Mella and Pasotti established the following two results relating extended irregular
dominating sets in cycles to strong starters in cyclic groups.

\begin{Lemma}{\rm \cite[Proposition~4.5]{MP}}\label{Domintion1}
Let $v\ge3$ be odd. Then $C_v$ admits an optimal extended irregular dominating set if and only if there is a strong starter in $\Z_v$.
\end{Lemma}

\begin{Lemma}{\rm \cite[Proposition~4.8]{MP}}\label{Domintion2}
Let $v\ge3$ be odd. If $C_v$ admits an optimal extended irregular dominating set, then so does $C_{2v}$.
\end{Lemma}

Combining Theorem \ref{abelian} with Lemmas \ref{Domintion1} and \ref{Domintion2}, and using explicit constructions for $C_6$, $C_{10}$, and $C_{18}$, we determine the extended irregular domination numbers of $C_v$ and $C_{2v}$ for all odd $v\ge3$, thereby resolving the existence problem posed in \cite{MP}.

\begin{Theorem}\label{Domintion-SS}
For every odd integer $v\ge3$, the cycle $C_v$ admits an optimal extended irregular dominating set if and only if $v\notin\{3,5,9\}$, in which case $\gamma_e(C_v)=(v+1)/2$. Moreover, $C_{2v}$ admits an optimal extended irregular dominating set for every odd integer $v\ge3$, with $\gamma_e(C_{2v})=v+1$.
\end{Theorem}

\begin{proof}
By Lemma \ref{Domintion1} and Theorem \ref{abelian}, the cycle $C_v$ admits an optimal extended irregular dominating set if and only if $v\notin\{3,5,9\}$. For every odd $v\ge3$ with $v\notin\{3,5,9\}$, Lemma \ref{Domintion2} guarantees that $C_{2v}$ admits an
optimal extended irregular dominating set. It remains to consider $C_{2v}$ for $v\in\{3,5,9\}$. For $C_6$ and $C_{10}$, Mella and Pasotti~\cite{MP} gave the following optimal labelings:
$$
C_6: (3,1,-,2,0,-),\qquad
C_{10}: (1,4,2,-,-,3,0,-,5,-),
$$
where the entries are listed in cyclic order and $-$ denotes an unlabeled vertex. For $C_{18}$, assign label $\ell$ to vertex $s_\ell\in\Z_{18}$ for $0\le\ell\le9$ as follows:
$$
\begin{array}{c|cccccccccc}
\ell&0&1&2&3&4&5&6&7&8&9\\ \hline
s_\ell&0&5&13&6&12&7&1&3&9&14
\end{array}
$$
The vertices $s_0,s_1,\ldots,s_9$ are pairwise distinct. For each $\ell$, the vertices covered by $s_\ell$ are, respectively,
$$
\{0\},\ \{4,6\},\ \{11,15\},\ \{3,9\},\
\{8,16\},\ \{2,12\},\ \{7,13\},\
\{10,14\},\ \{1,17\},\ \{5\}.
$$
These sets partition $\Z_{18}$, so this labeling gives an extended irregular dominating set of $C_{18}$.

We have thus established the existence of extended irregular dominating sets in $C_v$ for every odd $v\ge3$ with $v\notin\{3,5,9\}$, and in $C_{2v}$ for every odd $v\ge3$. Since cycles are vertex-transitive, \cite[Corollary~2.6]{MP} implies that every such set is
optimal and that its cardinality is one more than the diameter of the corresponding cycle. Consequently,
$$
\gamma_e(C_v)=\operatorname{diam}(C_v)+1
=\frac{v+1}{2}
$$
for every odd $v\ge3$ with $v\notin\{3,5,9\}$, and
$$
\gamma_e(C_{2v})=\operatorname{diam}(C_{2v})+1=v+1
$$
for every odd $v\ge3$.
\end{proof}

\subsection{Orthogonal starters and cyclic Room squares}

Strong starters can be used to construct pairwise orthogonal starters, which in turn provide a useful tool for constructing
Room squares and related combinatorial designs \cite{DS1992,ksv}.

Let $(G,+)$ be a finite abelian group of odd order with the identity $0_G$. Let $\mathcal S=\{\{s_i,t_i\}:1\le i\le (|G|-1)/2\}$ and $\mathcal T=\{\{u_i,v_i\}:1\le i\le (|G|-1)/2\}$ be two starters in $G$, indexed and oriented so that $s_i-t_i=u_i-v_i$ for $1\le i\le (|G|-1)/2$. Then $\mathcal S$ and $\mathcal T$ are said to be \emph{orthogonal} if
$$
s_i\ne u_i\quad(1\le i\le (|G|-1)/2)
\qquad\text{and}\qquad
s_j-u_j\ne s_i-u_i\quad(1\le i<j\le (|G|-1)/2).
$$
For example, $\mathcal S=\{\{2,3\},\{4,6\},\{1,5\}\}$ and $\mathcal T=\{\{3,4\},\{6,1\},\{5,2\}\}$ form a pair of orthogonal starters in $\Z_7$. 

For a starter $\mathcal S$ in $G$, define $-\mathcal S=\{\{-s,-t\}:\{s,t\}\in\mathcal S\}$. Also, define
\[
\mathcal P=\{\{x,-x\}:x\in G\setminus\{0_G\}\}.
\]
Both $-\mathcal S$ and $\mathcal P$ are starters in $G$.

\begin{Lemma}{\rm \cite[Theorem 2.1]{H}}\label{3-starter}
If $\mathcal S$ is a strong starter in $G$, then $\mathcal S$, $-\mathcal S$, and $\mathcal P$ are pairwise orthogonal starters.
\end{Lemma}
 
Applying Theorem \ref{abelian} and Lemma \ref{3-starter}, together with the known nonexistence results for the exceptional groups, gives the following characterization.

\begin{Theorem}\label{orthogonal-starters}
Let $G$ be a nontrivial finite abelian group of odd order. Then the following statements hold:
\begin{enumerate}
\item[$(1)$] There exists a pair of orthogonal starters in $G$ if and only if $G$ is not isomorphic to any of $\Z_3$, $\Z_5$, and $\Z_3\times\Z_3$.
\item[$(2)$] There exist three pairwise orthogonal starters in $G$ if and only if $G$ is not isomorphic to any of $\Z_3$, $\Z_5$, $\Z_9$, and $\Z_3\times\Z_3$.
\end{enumerate}
\end{Theorem}

\begin{proof}
Horton \cite[Section 4]{H} showed that $\Z_3$, $\Z_5$, and $\Z_3\times\Z_3$ admit no pair of orthogonal starters, whereas $\Z_9$ admits two orthogonal starters but no three pairwise orthogonal starters. For every other finite abelian group of odd order, Theorem~\ref{abelian} provides a strong starter, and Lemma \ref{3-starter} yields three pairwise orthogonal starters.
\end{proof}

We now apply Theorem \ref{orthogonal-starters} to cyclic Room squares. Let $X$ be a set of $v+1$ symbols. A \emph{Room square
of side $v$} is a $v\times v$ array $M$ whose cells are either empty or contain an unordered pair of
distinct symbols from $X$, such that each symbol occurs exactly once in every row and column of $M$, and every unordered
pair of distinct symbols from $X$ occurs in exactly one cell of $M$. Mullin and Wallis \cite{MW} showed that a Room square of side $v$ exists if and only if $v$ is odd and $v\not\in\{3,5\}$.

Taking $X=\Z_v\cup\{\infty\}$ and indexing the rows and columns of a Room square $M$ by $\Z_v$, we say that $M$ is \emph{cyclic} if
$$
M_{ii}=\{\infty,i\}\qquad(i\in\Z_v)
$$
and
$$
M_{i+1,j+1}=\{x+1:x\in M_{ij}\}
\qquad(i,j\in\Z_v),
$$
where the arithmetic is performed modulo $v$ and $\infty+1=\infty$.

\begin{Lemma} {\rm {\cite[Theorem~50.12]{handbookRoom}}}\label{pair}
There exists a pair of orthogonal starters in $\Z_v$ if and only if there exists a cyclic Room square of side $v$.
\end{Lemma}

Combining Lemma~\ref{pair} with Theorem~\ref{orthogonal-starters} completely determines the existence spectrum of cyclic Room squares.

\begin{Theorem}\label{Room square-cyclic}
For every integer $v>1$, a cyclic Room square of side $v$ exists if and only if $v$ is odd and $v\ge7$.
\end{Theorem}
\begin{proof}
Since the nonempty cells in each row of a Room square partition the $v+1$ symbols into pairs, $v$ must be odd. The conclusion then follows from Lemma~\ref{pair} and Theorem~\ref{orthogonal-starters}(1).
\end{proof}

\section*{Conflict of interest}

The authors declare no competing interests.

\section*{Data availability}

No datasets were generated or analysed during the current study.

\end{document}